\documentclass[11pt,leqno]{article}
\usepackage{amsmath}
\usepackage{amssymb}
\usepackage{dsfont}
\usepackage{multicol}
\usepackage{latexsym}
\usepackage[usenames]{color}
\usepackage[colorlinks=true,
linkcolor=webgreen,
filecolor=webbrown,
citecolor=webgreen]{hyperref}
\definecolor{webgreen}{rgb}{0,.5,0}
\definecolor{webbrown}{rgb}{.6,0,0}

\usepackage{graphics,amsmath,amssymb}
\usepackage{amsthm}
\usepackage{latexsym}

\def\N{{\mathds{N}}}
\def\Z{{\mathds{Z}}}

\newtheorem{thm}{Theorem}

\newtheorem{rem}{Remark}

\begin{document}

\title{\bf Some remarks on Ramanujan sums and cyclotomic polynomials}
\author{{\sc L\'aszl\'o T\'oth} \\ Department of Mathematics, University of P\'ecs \\
7624 P\'ecs, Ifj\'us\'ag u. 6., Hungary\\ \texttt{ltoth@gamma.ttk.pte.hu}}
\date{}
\maketitle

\begin{abstract} We investigate the polynomials
$\sum_{k=0}^{n-1} c_n(k)x^k$ and $\sum_{k=0}^{n-1} |c_n(k)| x^k$,
where $c_n(k)$ denote the Ramanujan sums. We point out connections
and analogies to the cyclotomic polynomials.
\end{abstract}




\section{Introduction}

The Ramanujan sum $c_n(k)$ is defined as the sum of $k$th powers
($k\in \Z$) of the primitive $n$th roots of unity, that is,
\begin{equation} \label{Ramanujan_sum}
c_n(k):=\sum_{j\in A_n} \eta_j^k,
\end{equation}
where $\eta_j=\exp(2\pi ij/n)$ and $A_n=\{j\in \N: 1\le j\le n,
(j,n)=1\}$. Here $c_n(k)$ is an $n$-periodic function of $k$, i.e.,
$c_n(k)=c_n(\ell)$ for any $k\equiv \ell$ (mod $n$). Note that for
$n\mid k$, $c_n(k)=c_n(0)=\varphi(n)$ is Euler's function and for
$(k,n)=1$, $c_n(k)=c_n(1)=\mu(n)$ is the M\"obius function.

The $n$th cyclotomic polynomial $\Phi_n(x)$ is the monic polynomial whose roots are the primitive
$n$th roots of unity, i.e.,
\begin{equation} \label{cyclotomic_pol}
\Phi_n(x):= \prod_{j\in A_n} \left(x- \eta_j \right).
\end{equation}

The following representations are well known:
\begin{equation} \label{Ramanujan_sum_Mobius}
c_n(k)= \sum_{d\mid (n,k)} d\mu(n/d),
\end{equation}
\begin{equation} \label{cyclotomic_pol_Mobius}
\Phi_n(x)= \prod_{d\mid n} (x^d-1)^{\mu(n/d)}.
\end{equation}

Cyclotomic polynomials and Ramanujan sums are closely related as it
is shown by the following theorem.

\begin{thm} \label{Th1} i) For any $n\ge 1$,
\begin{equation} \label{cyclotomic_Ramanujan_deriv}
(x^n-1)\frac{\Phi'_n(x)}{\Phi_n(x)} = \sum_{k=1}^n c_n(k)x^{k-1},
\end{equation}
where $\Phi'_n(x)$ is the derivative of $\Phi_n(x)$.

ii) For any $n>1$ and $|x|<1$,
\begin{equation} \label{cyclotomic_Ramanujan_exp}
\Phi_n(x) = \exp \left(-\sum_{k=1}^{\infty} \frac{c_n(k)}{k}x^k
\right).
\end{equation}
\end{thm}

These formulae are not widely known and were first derived by Nicol
\cite[Th.\ 3.1, Cor.\ 3.2]{Nic1962}. In that paper formula
\eqref{cyclotomic_Ramanujan_deriv} was deduced by differentiating
\eqref{cyclotomic_pol_Mobius}, which gives
\begin{equation} \label{cyclotomic_sum}
\frac{\Phi'_n(x)}{\Phi_n(x)} = \sum_{d\mid n} \frac{d\mu(n/d)x^{d-1}}{x^d-1} \quad (n\ge 1)
\end{equation}
and then using \eqref{Ramanujan_sum_Mobius}, while
\eqref{cyclotomic_Ramanujan_exp} was obtained as a corollary of
\eqref{cyclotomic_Ramanujan_deriv}.

Formula \eqref{cyclotomic_sum} was given also by Motose \cite[Th.\
1]{Mot2005}, without referring to the paper of Nicol \cite{Nic1962}.
The following result was  obtained by Motose \cite[Lemma 1, Th.\
1]{Mot2005} in the same paper.

\begin{thm} \label{Th2} For any $n\ge 1$,
\begin{equation} \label{cyclotomic_x_inverse}
\frac{\Phi'_n(1/x)}{x\Phi_n(1/x)}= \frac1{1-x^n} \sum_{k=0}^{n-1}
c_n(k)x^k = \sum_{d\mid n} \frac{d\mu(n/d)}{1-x^d}.
\end{equation}
\end{thm}

Note that \eqref{cyclotomic_x_inverse} is a simple consequence of
formulae \eqref{cyclotomic_Ramanujan_deriv} and
\eqref{cyclotomic_sum} by putting $x:=1/x$.

In this paper we first give new direct proofs of
\eqref{cyclotomic_Ramanujan_deriv} and
\eqref{cyclotomic_Ramanujan_exp} which use only the definitions of
the Ramanujan sums and of the cyclotomic polynomials (Section 2).

Then in Section 3 we investigate the polynomials with integer
coefficients
\begin{equation} \label{R_n(x)}
R_n(x):=\sum_{k=0}^{n-1} c_n(k)x^k
\end{equation}
appearing in \eqref{cyclotomic_x_inverse}. We deduce for $R_n(x)$
formulas which are similar to the following well known formulas
valid for the cyclotomic polynomials:
$\Phi_n(x)=\Phi_{\gamma(n)}(x^{n/\gamma(n)})$, where
$\gamma(n)=\prod_{p\mid n} p$ is the squarefree kernel of $n$,
$\Phi_{np}(x)=\Phi_{n}(x^p)$ for any prime $p\mid n$,
$\Phi_{np}(x)=\Phi_{n}(x^p)/\Phi_n(x)$ for any prime $p\nmid n$,
$\Phi_{2n}(x)=\Phi_{n}(-x)$ for any $n$ odd, see for ex.
\cite{Tha2000}. We also derive certain divisibility properties of
the polynomials $R_n(x)$.

In Section 4 we consider the polynomials
\begin{equation} \label{T_n(x)}
T_n(x):=\sum_{k=0}^{n-1} |c_n(k)| x^k
\end{equation}
and compare their properties to those of the polynomials $R_n(x)$.

We show -- among others -- that $R_{2n}(x)=(1-x^n)R_n(-x)$,
$T_{2n}(x)=(1+x^n)T_n(x)$ for any $n\ge 1$ odd, and that the
cyclotomic polynomial $\Phi_n(x)$ divides $T_n(x)$ for any $n\ge 2$
even.

Section 5 contains tables of the polynomials $R_n(x)$ and $T_n(x)$
for $1\le n\le 20$.

For material concerning Ramanujan sums we refer to the books
\cite{HarWri1975,McC1986,MonVau2007}.


\section{Proof of Theorem 1}

\begin{proof} i) First note that for the generating function of the sequence
$(c_n(k))_{k\ge 1}$ we have by using the periodicity of the
Ramanujan sums,
\begin{equation*}
\sum_{k=1}^{\infty} c_n(k)x^k = \sum_{\ell=0}^{\infty} \sum_{j=1}^n
c_n(\ell n+j)x^{\ell n+j} = \sum_{\ell=0}^{\infty} x^{\ell n}
\sum_{j=1}^n c_n(j) x^{j} =\frac1{1-x^n} \sum_{j=1}^n c_n(j)x^j.
\end{equation*}

Now let $|x|<1$. Applying the power series
$(1-t)^{-1}=1+t+t^2+\ldots$ for $t=x/\eta_j$, where $|t|=|x|<1$,
\begin{equation*}
\frac{\Phi'_n(x)}{\Phi_n(x)} = \sum_{j\in A_n} \frac1{x-\eta_j} = -
\sum_{j\in A_n} \frac1{\eta_j}\cdot \frac1{1-x/\eta_j} = -
\sum_{j\in A_n} \frac1{\eta_j} \sum_{k=0}^{\infty}
\frac{x^k}{\eta_j^k}
\end{equation*}
\begin{equation*}
= - \sum_{k=0}^{\infty} x^k  \sum_{j\in A_n} \eta_j^{-k-1} = -
\sum_{k=0}^{\infty} x^k  c_n(-k-1)= - \sum_{k=1}^{\infty}
c_n(-k)x^{k-1}
\end{equation*}
\begin{equation*}
= - \sum_{k=1}^{\infty} c_n(k) x^{k-1} = \frac1{x^n-1} \sum_{j=1}^n
c_n(j)x^{j-1},
\end{equation*}
where we have used that $c_n(-k)=c_n(k)$ for any $k$. Justification:
if in \eqref{Ramanujan_sum} $j$ runs through a reduced residue
system (mod $n$), then so does $-j$. Hence the given polynomial
identity holds, which finishes the proof of i).

ii) We use that for $n>1$,
\begin{equation} \label{cyclotomic_pol_var}
\Phi_n(x)=\prod_{j\in A_n} \left(1-\frac{x}{\eta_j} \right).
\end{equation}

This follows from \eqref{cyclotomic_pol} by $\prod_{j\in A_n}
\eta_j=1$, valid for $n>2$. Note that \eqref{cyclotomic_pol_var}
holds also for $n=2$. We have, using the power series
$\log(1-t)=-t-t^2/2-t^3/3-\ldots$ for $t=x/\eta_j$, where
$|t|=|x|<1$,
\begin{equation*}
\log \Phi_n(x)= \sum_{j\in A_n} \log \left(1-\frac{x}{\eta_j}
\right) = - \sum_{j\in A_n} \sum_{k=1}^{\infty}
\frac{x^k}{k\eta_j^k}= - \sum_{k=1}^{\infty} \frac{x^k}{k}
\sum_{j\in A_n} \eta_j^{-k}
\end{equation*}
\begin{equation*}
=  - \sum_{k=1}^{\infty} \frac{x^k}{k} c_n(-k) = -
\sum_{k=1}^{\infty} \frac{x^k}{k} c_n(k).
\end{equation*}

Alternatively, one can apply that for $n>1$,
\begin{equation}
\label{cyclotomic_pol_Mobius_var} \Phi_n(x)= \prod_{d\mid n}
(1-x^d)^{\mu(n/d)},
\end{equation}
which follows at once by \eqref{cyclotomic_pol_Mobius} and by
$\sum_{d\mid n} \mu(n/d)=0$ ($n>1$). We deduce
\begin{equation*}
\log \Phi_n(x)= \sum_{d\mid n} \mu(n/d) \log (1-x^d) = - \sum_{d\mid
n} \mu(n/d) \sum_{j=1}^{\infty} \frac{x^{dj}}{j}
\end{equation*}
\begin{equation*}
= - \sum_{m=1}^{\infty} \frac{x^m}{m}\sum_{d\mid (n,m)} d \mu(n/d) =
-\sum_{m=1}^{\infty} \frac{x^m}{m}c_n(m),
\end{equation*}
using \eqref{Ramanujan_sum_Mobius}. This approach was given and
applied by Erd\H os and Vaughan \cite[Proof of Th.\ 1]{ErdVau1974}.
\end{proof}


\section{The polynomials $R_n(x)$}

In this section we investigate properties of the polynomials
$R_n(x)$ defined by \eqref{R_n(x)}. Note that the polynomials
appearing in \eqref{cyclotomic_Ramanujan_deriv} are given by
$P_n(x):=\sum_{k=1}^n c_n(k)x^{k-1}$. The connection between the
polynomials $R_n(x)$ and $P_n(x)$ is given by
$xP_n(x)=R_n(x)+\varphi(n)(x^n-1)$. Hence it is sufficient to study
the polynomials $R_n(x)$.

According to \eqref{cyclotomic_x_inverse} for any $n\ge 1$,
\begin{equation} \label{R_n(x)_representation}
R_n(x) = (1-x^n) \sum_{d\mid n} \frac{d\mu(n/d)}{1-x^d}.
\end{equation}

\begin{thm} \label{Th3} Let $n\ge 1$.

i) The number of nonzero coefficients of $R_n(x)$ is $\gamma(n)$.

ii) The degree of $R_n(x)$ is $n-n/\gamma(n)$.

iii) $R_n(x)$ has coefficients $\pm 1$ if and only if $n$ is
squarefree and in this case the number of coefficients $\pm 1$ of
$R_n(x)$ is $\varphi(n)$ for $n$ odd and is $2\varphi(n/2)$ for $n$
even.
\end{thm}

\begin{proof} For $n=1$ the assertions hold true.
Let $n=p_1^{a_1}\cdots p_r^{a_r}>1$.

i) We use that $c_n(k)$ is multiplicative in $n$ and for any prime
power $p^a$,
\begin{equation*} \label{c_n(k)}
c_{p^a}(k)=\begin{cases} p^a-p^{a-1}, \ \text{ if } \ p^a\mid k, \\
-p^{a-1}, \ \text{ if } \ p^{a-1}\mid k, p^a\nmid k, \\ 0, \ \text{
if } \ p^{a-1}\nmid k.
\end{cases}
\end{equation*}

Therefore, $c_n(k)\ne 0$ if and only if $p_1^{a_1-1}\mid k$, ...,
$p_r^{a_r-1}\mid k$, i.e., $k=p_1^{a_1-1}\cdots p_r^{a_r-1}m$ with
$0\le m < p_1\cdots p_r= \gamma(n)$. Hence the number of nonzero
values of $c_n(k)$ is $\gamma(n)$.

ii) By the proof of i) the largest $k$ such that $c_n(k)\ne 0$ is
$k=p_1^{a_1-1}\cdots p_r^{a_r-1}(p_1\cdots p_r-1)=n-n/\gamma(n)$,
and this is the degree of $R_n(x)$.

iii) $c_n(k)=\pm 1$ if and only if $c_{p_i^{a_i}}(k)=\pm 1$ for any
$i\in \{1,\ldots, r\}$, that is $a_i=1$ for any $i$ ($n$ is
squarefree) and either $p_i\nmid k$ or $p_i=2\mid k$ for any $i$.

Suppose that $n=p_1\cdots p_r$ (squarefree). If $n$ is odd, then by
condition $p_i\nmid k$ for any $i$ we have $(n,k)=1$, hence the
number of such values of $k$ is $\varphi(n)$. For $n$ even either
$(k,n)=1$ or $k=2\ell$ with $(\ell,n/2)=1$. We obtain that the
number of such values of $k$ is $\varphi(n)+ \varphi(n/2)=
2\varphi(n/2)$.
\end{proof}

We have for any $n>1$, $R_n(0)=c_n(0)=\varphi(n)$ and
$R_n(1)=\sum_{k=0}^{n-1} c_n(k)=0$, as it is well known. Hence $1-x$
divides $R_n(x)$ for any $n>1$. Now a look at the polynomials
$R_n(x)$, see Section 5, suggests that $1+x$ divides $R_n(x)$ for
any $n>2$ even. This is confirmed by the next result.

\begin{thm} \label{Th4} We have $R_2(-1)=2$ and

i) $R_n(-1)=\varphi(n)$ for any $n\ge 1$ odd,

ii) $R_n(-1)=0$ for any $n>2$ even,

iii) the cyclotomic polynomial $\Phi_n(x)$ divides the polynomial
$R_n(x)-n$ for any $n\ge 1$.
\end{thm}

\begin{proof} i) Consider also the polynomials
$Q_n(x)=\sum_{k=0}^n c_n(k)x^k =R_n(x)+\varphi(n)x^n$, which are
symmetric for any $n\ge 1$ since $c_n(k)=c_n(n-k)$ ($0\le k\le n$).
Hence, for any $n$ odd, $Q_n(-1)=0$ and $R_n(-1)=\varphi(n)$.

ii) Now use \eqref{R_n(x)_representation}, which can be written as
\begin{equation}
R_n(x)= \sum_{d\mid n} d\mu(n/d)(x^{n-d}+x^{n-2d}+\ldots +x^d+1).
\end{equation}

We obtain that for any $n=2k>2$ even,
\begin{equation*}
R_n(-1)= \sum_{\substack{d\mid n\\d \text{ even}}}
d\mu(n/d)\frac{n}{d}= n\sum_{\substack{d\mid n\\d \text{ even}}}
\mu(n/d)= n\sum_{\delta \mid k} \mu(k/\delta)=0.
\end{equation*}

iii) If $\eta$ is any primitive $n$th root of unity, then
$R_n(\eta)=n$. This follows from \eqref{R_n(x)_representation}:
\begin{equation*}
R_n(x)= n+ (1-x^n) \sum_{\substack{d\mid n\\ d<n}}
\frac{d\mu(n/d)}{1-x^d},
\end{equation*}
where $\eta^d\ne 1$ for any $d\mid n$, $d<n$.
\end{proof}

If $p$ is a prime, then it follows from
\eqref{R_n(x)_representation} that
\begin{equation} \label{R_prim}
R_p(x)= (p-1)-x-x^2-\ldots -x^{p-1}.
\end{equation}

Also, if $p,q$ are distinct primes, then
\begin{equation} \label{R_pq}
R_{pq}(x)= (p-1)(q-1)+ x+x^2+\ldots + x^{pq-1}
\end{equation}
\begin{equation*}- p(x^p+x^{2p}+\ldots + x^{(q-1)p})- q(x^q+x^{2q}+\ldots +x^{(p-1)q}).
\end{equation*}

Next we show that $R_n(x)$ have some properties which are similar to
those of the cyclotomic polynomials $\Phi_n(x)$.

\begin{thm} \label{Th5} i) If $n\ge 1$, then
\begin{equation} \label{R_kernel_representation} R_n(x)=
\frac{n}{\gamma(n)} R_{\gamma(n)} (x^{n/\gamma(n)}).
\end{equation}

ii) Let $n\ge 1$ and $p$ be a prime. If $p\mid n$, then
$R_{np}(x)=pR_n(x^p)$. If $p\nmid n$, then
\begin{equation} \label{R_np}
R_{np}(x)= pR_n(x^p) - (1+x^n+x^{2n}+\ldots +x^{(p-1)n})R_n(x).
\end{equation}

iii) If $n>1$, $p$ is a prime and $p\nmid n$, then $(1-x^p) \mid
R_{np}(x)$.
\end{thm}

\begin{proof} i) From
\eqref{R_n(x)_representation} again,
\begin{equation*}
R_n(x)= (1-x^n) \sum_{d\mid n} \frac{\mu(d)\frac{n}{d}}{1-x^{n/d}} =
(1-x^n) \sum_{d\mid \gamma(n)} \frac{\mu(d)\frac{n}{d}}{1-x^{n/d}},
\end{equation*}
and by this representation of $R_n(x)$,
\begin{equation*}
R_{\gamma(n)} (x^{n/\gamma(n)}) =(1-x^n) \sum_{d\mid \gamma(n)}
\frac{\mu(d)\frac{\gamma(n)}{d}}{1-x^{n/d}}=
(1-x^n)\frac{\gamma(n)}{n} \sum_{d\mid \gamma(n)}
\frac{\mu(d)\frac{n}{d}}{1-x^{n/d}}= \frac{\gamma(n)}{n}R_n(x).
\end{equation*}

ii) For $p\mid n$ this follows at once from i) by
$\gamma(np)=\gamma(n)$. Now let $p\nmid n$. Then by
\eqref{R_n(x)_representation},
\begin{equation*}
R_{np}(x) = (1-x^{np}) \sum_{d\mid np} \frac{d\mu(np/d)}{1-x^d}=
(1-x^{np})\left( \sum_{d\mid n} \frac{d\mu(np/d)}{1-x^d}+
\sum_{d=\delta p, \ \delta \mid n} \frac{\delta p
\mu(n/\delta)}{1-x^{\delta p}}\right)
\end{equation*}
\begin{equation*}
= (1-x^{np})\left( - \sum_{d\mid n} \frac{d\mu(n/d)}{1-x^d} + p
\sum_{\delta \mid n} \frac{\delta \mu(n/\delta)}{1-x^{\delta
p}}\right) = p R_n(x^p) - \frac{1-x^{np}}{1-x^n} R_n(x). \ \Box
\end{equation*}

iii) Using that $x=1$ is a root of $R_n(x)$ for $n>1$ we deduce that
$(1-x^p)\mid R_n(x^p)$ and by the formula \eqref{R_np} we obtain
$(1-x^p)\mid R_{np}(x)$.
\end{proof}

In particular, for any prime power $p^k$ ($k\ge 1$),
\begin{equation} \label{R_p^k}
R_{p^k}(x) = p^{k-1} R_p(x^{p^{k-1}})=
p^{k-1}(p-1-x^{p^{k-1}}-x^{2p^{k-1}}-\ldots -x^{(p-1)p^{k-1}})
\end{equation}
and for $p=2$,
\begin{equation} \label{eq_2^k}
R_{2^k}(x)=2^{k-1}(1-x^{2^{k-1}}).
\end{equation}

\begin{thm} \label{Th6} i) For any $n\ge 1$ odd,
\begin{equation} \label{n_odd}
R_{2n}(x)=(1-x^n) R_n(-x),
\end{equation}

ii) More generally, for any $n\ge 1$ odd and any $k\ge 1$,
\begin{equation} \label{2^k}
R_{2^k n}(x)=2^{k-1} (1-x^{2^{k-1}n}) R_n(-x^{2^{k-1}}).
\end{equation}
\end{thm}

\begin{proof} i) By \eqref{R_np} we have
\begin{equation*}
R_{2n}(x)=2R_n(x^2)-(1+x^n)R_n(x)
\end{equation*}
\begin{equation*}
= 2(1-x^{2n}) \sum_{d\mid n} \frac{d\mu(n/d)}{1-x^{2d}}-
(1+x^n)(1-x^n) \sum_{d\mid n} \frac{d\mu(n/d)}{1-x^d}
\end{equation*}
\begin{equation*}
= (1-x^{2n}) \sum_{d\mid n} \frac{d\mu(n/d)}{1+x^d}=(1-x^n)R_n(-x),
\end{equation*}
hence \eqref{n_odd} holds.

ii) By \eqref{R_kernel_representation} and \eqref{n_odd} we obtain
\begin{equation*}
R_{2^k n}(x)=
\frac{2^{k-1}n}{\gamma(n)}R_{2\gamma(n)}(x^{2^{k-1}n/\gamma(n)})=
\frac{2^{k-1}n}{\gamma(n)}(1-x^{2^{k-1}n})
R_{\gamma(n)}(-x^{2^{k-1}n/\gamma(n)}).
\end{equation*}

Here by \eqref{R_kernel_representation} again,
\begin{equation*}
R_n(-x^{2^{k-1}})= \frac{n}{\gamma(n)}
R_{\gamma(n)}(-x^{2^{k-1}n/\gamma(n)}),
\end{equation*}
ending the proof.
\end{proof}

\begin{thm} \label{Th7}

i) If $n=p^k$, $p$ prime, $k\ge 1$, then
\begin{equation} \label{div_1}
(1-x^{p^{k-1}}) \mid R_n(x).
\end{equation}

ii) If $n=2^k m$, $k\ge 1$, $m>1$ odd, then
\begin{equation} \label{div_2}
(1-x^{n/2})(1+x^{n/\gamma(n)}) \mid R_n(x).
\end{equation}

iii) If $n=p^k m$, $p>2$ prime, $k\ge 1$, $m>1$ odd, $p\nmid m$,
then
\begin{equation} \label{div_3}
(1-x^{pn/\gamma(n)}) \mid R_n(x).
\end{equation}

iv) If $n=2^k m$, $k\ge 1$, $m>1$ odd, $m$ has at least two prime
divisors, $p$ prime, $p\mid m$, then
\begin{equation} \label{div_4}
(1-x^{n/2})(1+x^{pn/\gamma(n)}) \mid R_n(x).
\end{equation}
\end{thm}

\begin{proof} i) $R_{p^k}(x)=p^{k-1}R_p(x^{p^{k-1}})$ by \eqref{R_p^k},
and use that $x=1$ is a root of $R_p(x)$.

ii) For $n=2^k m$, $k\ge 1$, $m>1$ odd we have
\begin{equation*}
R_n(x) =\frac{2^{k-1}m}{\gamma(m)} \left(1-x^{2^{k-1}m}\right)
R_{\gamma(m)}(-x^{2^{k-1}m/\gamma(m)}),
\end{equation*}
see the proof of Theorem 6/ii), and use that $x=1$ is a root of
$R_{\gamma(m)}(x)$.

iii) Now
\begin{equation*}
R_n(x) =\frac{n}{\gamma(n)} R_{p\gamma(m)}(x^{p^{k-1}m/\gamma(m)}),
\end{equation*}
where $(1-x^p) \mid R_{p\gamma(m)}(x)$ with $\gamma(m)>1$, cf.
Theorem 5/iii).

iv) By combining the above results.
\end{proof}

As examples, Theorem 7 gives that $(1-x^9)(1+x^3)\mid R_{18}(x)$ and
$(1-x^{15})(1+x^3)\mid R_{30}(x)$. It is possible to deduce from
Theorem 7 other divisibility properties for the polynomials
$R_n(x)$, e.g. the next one.

\begin{thm} \label{Th8} For any $k\ge 1$ and $m>1$,
\begin{equation} \label{div_5}
\left(1+x^{2^{k-1}}\right) \mid R_{2^k m}(x).
\end{equation}
\end{thm}

\begin{proof} Follows from Theorem 7/ii).
\end{proof}

Another representation of the polynomials $R_n(x)$ is given by

\begin{thm} \label{Th9} For any $n\ge 1$,
\begin{equation} \label{P(x)_representation}
R_n(x)= \varphi(n) \left(1-x^n+ \sum_{d\mid n}
\frac{\mu(d)}{\varphi(d)} \Psi_d(x^{n/d})\right),
\end{equation}
where $\Psi_n(x)=\sum_{j\in A_n} x^j$.
\end{thm}

\begin{proof} We use H\"older's formula
\begin{equation}
c_n(k)=\frac{\varphi(n)\mu(n/(n,k))}{\varphi(n/(n,k))}
\end{equation}
and by grouping the terms according to $(n,k)=d$, obtain
\begin{equation*}
\sum_{k=1}^n c_n(k)x^k = \sum_{k=1}^n
\frac{\varphi(n)\mu(n/(n,k))}{\varphi(n/(n,k))} x^k = \varphi(n)
\sum_{d\mid n} \frac{\mu(n/d)}{\varphi(n/d)} \sum_{j \in A_{n/d}}
x^{dj}
\end{equation*}
\begin{equation*}
= \varphi(n) \sum_{d\mid n} \frac{\mu(d)}{\varphi(d)} \Psi_d
(x^{n/d}).
\end{equation*}
\end{proof}

\begin{rem} {\rm The polynomials $\Psi_n(x)=\sum_{j\in A_n} x^j$
are $\Psi_1(x)=x$, $\Psi_2(x)=x$, $\Psi_3(x)=x+x^2$,
$\Psi_4(x)=x+x^3$, $\Psi_5(x)=x+x^2+x^3+x^4$, $\Psi_6(x)=x+x^5$,
etc., having the representations
\begin{equation} \label{Psi}
\Psi_n(x)= (1-x^n) \sum_{d\mid n}  \frac{\mu(d)x^d}{1-x^d}= (1-x^n)
\sum_{d\mid n}  \frac{\mu(d)}{1-x^d}, \quad n>1,
\end{equation}
the first one being valid for $n\ge 1$.

If $\eta$ is any primitive $n$th root of unity, then
$\Psi_n(\eta)=\mu(n)$. Hence the cyclotomic polynomial $\Phi_n(x)$
divides the polynomial $\Psi_n(x)-\mu(n)$ for any $n\ge 1$. For
these properties see \cite[p.\ 71]{Wil1994}. Furthermore, it is
immediate from \eqref{Psi} that $\Psi_n(1)=\varphi(n)$ for any $n\ge
1$ and $\Psi_n(-1)= - \varphi(n)$ for any $n\ge 2$ even. Also,
$\Psi_1(-1)=-1$ and $\Psi_n(-1)=0$ for any $n>1$ odd, since by
\eqref{Psi}, $\Psi_n(-1)= (1-(-1)^n) \sum_{d\mid n}
\frac{\mu(d)}{1-(-1)^d}= - \sum_{d\mid n} \mu(d)=0$.}
\end{rem}


\section{The polynomials $T_n(x)$}

We consider in what follows the polynomials $T_n(x)$ given by
\eqref{T_n(x)}.

Theorem 3 holds for the polynomials $T_n(x)$, as well. Also, for any
prime $p$, $T_p(x)=p-1+x+x^2+\ldots +x^{p-1}$, which follows at once
from \eqref{R_prim}.

\begin{thm} \label{Th10}
\begin{equation} \label{T_n(x)_representation}
T_n(x)= \varphi(n)\left(1-x^n + \sum_{d\mid n}
\frac{\mu^2(d)}{\varphi(d)} \Psi_d(x^{n/d})\right).
\end{equation}
\end{thm}

\begin{proof} Similar to the proof of Theorem 9,
using H\"older's formula.
\end{proof}

Note that for every $n\ge 1$, $T_n(0)=|c_n(0)|=\varphi(n)$ is
Euler's function. Also, $T_n(1)=\sum_{k=0}^{n-1}
|c_n(k)|=\varphi(n)2^{\omega(n)}$, where $\omega(n)$ denotes, as
usual, the number of distinct prime factors of $n$. This identity
follows at once by \eqref{T_n(x)_representation} and is given in
\cite{Bac1997}.

Now we deduce for $T_n(x)$ a formula which is similar to
\eqref{R_n(x)_representation}.

\begin{thm} \label{Th11} For any $n\ge 1$,
\begin{equation} \label{T_n_formula}
T_n(x)= (1-x^n)\varphi(n) \sum_{d\mid n}
\frac{\mu^2(d)f_d(n/d)}{\varphi(d)(1-x^{n/d})},
\end{equation}
where $f_k(n)$ denotes the multiplicative function in $n$ given by
\begin{equation}
f_k(n)=\prod_{\substack{p\mid n\\ p\nmid k}}
\left(1-\frac1{p-1}\right).
\end{equation}
\end{thm}

Note that $f_k(n)=0$ for any $n$ even and $k$ odd.

\begin{proof} Formula \eqref{Ramanujan_sum_Mobius} can
not be used in this case and we start with
\eqref{T_n(x)_representation}. Using also \eqref{Psi} we deduce
\begin{equation*}
T_n(x)- (1-x^n)\varphi(n) = \varphi(n) \sum_{d\mid n}
\frac{\mu^2(d)}{\varphi(d)} \left(1-(x^{n/d})^d\right) \sum_{\delta
\mid d} \frac{\mu(\delta)x^{n\delta/d}}{1-x^{n\delta/d}}
\end{equation*}
\begin{equation*}
=\varphi(n)(1-x^n) \sum_{ab\delta=n} \frac{\mu^2(b\delta)
\mu(\delta)x^{a\delta}}{\varphi(b\delta)(1-x^{a\delta})}
\end{equation*}
\begin{equation*}
=\varphi(n)(1-x^n) \sum_{\substack{ab\delta=n\\ (b,\delta)=1}}
\frac{\mu^2(b)\mu(\delta)
x^{a\delta}}{\varphi(b)\varphi(\delta)(1-x^{a\delta})}
\end{equation*}
\begin{equation*}
=(1-x^n)\varphi(n) \sum_{bc=n} \frac{\mu^2(b)x^c}{\varphi(b)(1-x^c)}
\sum_{\substack{a\delta=n\\(\delta,b)=1}}
\frac{\mu(\delta)}{\varphi(\delta)},
\end{equation*}
where the inner sum is $f_b(c)$ and the given formula follows by
writing $\frac{x^c}{1-x^c}= -1+\frac1{1-x^c}$ and using that
$\sum_{bc=n} \frac{\mu^2(b)f_b(c)}{\varphi(b)}=1$, which can be
checked easily by the multiplicativity of the involved functions.
\end{proof}

In particular, if $p,q$ are distinct primes, then by
\eqref{T_n_formula} we obtain
\begin{equation}
T_{pq}(x)= (p-1)(q-1)+ x+x^2+\ldots + x^{pq-1}
\end{equation}
\begin{equation*}+ (p-2)(x^p+x^{2p}+\ldots + x^{(q-1)p})+ (q-2)(x^q+x^{2q}+\ldots
+x^{(p-1)q}),
\end{equation*}
which follows also from \eqref{R_pq}.

\begin{thm} \label{Th12} We have

i) $T_n(-1)=\varphi(n)$ for any $n\ge 1$ odd,

ii) $T_n(-1)=0$ for any $n=4k+2$, $k\ge 0$,

iii) $T_n(-1)=\varphi(n)2^{\omega(n)}$ for any $n=4k$, $k\ge 1$,

iv) $T_n(\eta)= n \prod_{p\mid n} (1-\frac{2}{p})$ for any primitive
$n$th root of unity $\eta$. The cyclotomic polynomial $\Phi_n(x)$
divides the polynomial $T_n(x)$ for any $n\ge 2$ even.
\end{thm}

\begin{proof} For i)--iii) we use formula
\eqref{T_n(x)_representation} and the properties of the polynomials
$\Psi_n(x)$, mentioned in Remark 1.

i) For any $n\ge 1$ odd,
\begin{equation*}
T_n(-1)=\varphi(n) \left(2+\Psi_1(-1) + \varphi(n)\sum_{d\mid n, \,
d>1} \frac{\mu^2(d)}{\varphi(d)} \Psi_d(-1) \right) = \varphi(n).
\end{equation*}

ii) For any $n=4k+2$, $k\ge 0$,
\begin{equation*}
T_n(-1) = \varphi(n) \sum_{d\mid n} \frac{\mu^2(d)}{\varphi(d)}
\Psi_d((-1)^{n/d})
\end{equation*}
\begin{equation*}
=\varphi(n) \sum_{d\mid 2k+1} \frac{\mu^2(d)}{\varphi(d)} \Psi_d(1)
+ \varphi(n) \sum_{d=2\delta, \, \delta \mid 2k+1}
\frac{\mu^2(2\delta)}{\varphi(2\delta)} \Psi_{2\delta}(-1)
\end{equation*}
\begin{equation*} =  \varphi(n) \sum_{d\mid 2k+1} \mu^2(d)
-\varphi(n) \sum_{\delta \mid 2k+1} \mu^2(\delta) =0.
\end{equation*}

iii) For any $n=4k$, $k\ge 1$,
\begin{equation*}
T_n(-1) = \varphi(n) \sum_{d\mid n} \frac{\mu^2(d)}{\varphi(d)}
\Psi_d((-1)^{n/d}),
\end{equation*}
where for any $d$ with $4\mid d$, $\mu^2(d)=0$. Hence
\begin{equation*}
T_n(-1) = \varphi(n) \sum_{d\mid n, \, 4\nmid d}
\frac{\mu^2(d)}{\varphi(d)} \Psi_d(1)= \varphi(n) \sum_{d\mid n, \,
4\nmid d} {\mu^2(d)} = \varphi(n) \sum_{d\mid n} {\mu^2(d)}=
\varphi(n)2^{\omega(n)}.
\end{equation*}

iv) Now we use \eqref{T_n_formula}. The property follows from
\begin{equation*}
T_n(x) = \varphi(n)f_1(n) + (1-x^n) \sum_{d\mid n, d>1}
\frac{\mu^2(d)f_d(n/d)}{\varphi(d)(1-x^{n/d})},
\end{equation*}
where $\eta^k \ne 1$ for any $k\mid n$, $k<n$ and $\varphi(n)
f_1(n)=n\prod_{p\mid n} (1-2/p)$.
\end{proof}

\begin{thm} \label{Th13} i) If $n\ge 1$, then
\begin{equation} \label{kernel_representation}
T_n(x)= \frac{n}{\gamma(n)} T_{\gamma(n)} (x^{n/\gamma(n)}).
\end{equation}

ii) Let $n\ge 1$ and $p$ be a prime. If $p\mid n$, then
$T_{np}(x)=pT_n(x^p)$. If $p\nmid n$, then $\displaystyle T_{np}(x)=
(p-2)\varphi(n) T_n(x^p) + (1+x^n+x^{2n}+\ldots +x^{(p-1)n})
T_n(x)$.
\end{thm}

\begin{proof} i) This follows at once from Theorem 5/i) and
from the definitions of the polynomials $T_n(x)$ and $R_n(x)$.

ii) For $p\mid n$ this follows from i) by $\gamma(np)=\gamma(n)$.
For $p\nmid n$ by \eqref{T_n_formula},
\begin{equation*}
T_{np}(x)= (1-x^{np})\varphi(np) \sum_{d\mid np}
\frac{\mu^2(d)f_d(np/d)}{\varphi(d)(1-x^{np/d})}
\end{equation*}
\begin{equation*}
= (1-x^{np})\varphi(n)\varphi(p) \left( \sum_{d\mid n}
\frac{\mu^2(d)f_d(n/d)f_d(p)}{\varphi(d)(1-x^{np/d})}+
\sum_{d=\delta p, \ \delta \mid n} \frac{\mu^2(\delta p)f_{\delta
p}(n/\delta)}{\varphi(\delta p)(1-x^{n/\delta})}\right)
\end{equation*}
\begin{equation*}
= (1-x^{np})\varphi(n) \left((p-2) \sum_{d\mid n}
\frac{\mu^2(d)f_d(n/d)}{\varphi(d)(1-x^{np/d})}+ \sum_{\delta \mid
n} \frac{\mu^2(\delta)f_{\delta p}(n/\delta)}{\varphi(\delta)
(1-x^{n/\delta})}\right)
\end{equation*}
\begin{equation*}
= (p-2)\varphi(n) T_n(x^p) + \frac{1-x^{np}}{1-x^n}T_n(x),
\end{equation*}
where $f_{\delta p}(n/\delta)=f_{\delta}(n/d)$ for any $\delta \mid
n$, $p\nmid n$. \end{proof}

\begin{thm} \label{Th14}
i) For any $n\ge 1$ odd, $T_{2n}(x)=(1+x^n) T_n(x)$.

ii) For any $n\ge 1$ odd and any $k\ge 1$,
\begin{equation} \label{2^k_T}
T_{2^k n}(x)=2^{k-1}\left(1+x^{2^{k-1}n}\right) T_n(x^{2^{k-1}}).
\end{equation}

iii) For any even $n$, $\left(1+x^{n/2}\right) \mid T_n(x)$.
\end{thm}

\begin{proof} i) This follows at once from Theorem 13/ii) by $p=2$.

ii), iii) The same proof as for the polynomials $R_n(x)$.
\end{proof}

\begin{rem} {\rm Consider the polynomials
\begin{equation}
V_n(x)=\sum_{k=0}^{n-1} (c_n(k))^2 x^k.
\end{equation}

For every $n\ge 1$, $V_n(0)=(c_n(0))^2=(\varphi(n))^2$ and
$V_n(1)=\sum_{k=0}^{n-1} (c_n(k))^2 = n\varphi(n)$, as it is known.
For the polynomials $V_n(x)$ similar properties can be derived as
for $T_n(x)$.}
\end{rem}


\section{Tables of $R_n(x)$ and $T_n(x)$}

The next two tables were produced using Maple. The polynomials
$R_n(x)$ were generated by the following procedure (similar for
$T_n(x)$):
\begin{verbatim}
with(numtheory): Ramanujanpol:= proc(n,x) local a, k: a:= 0: for k
from 0 to n-1 do a:=a+phi(n)*mobius(n/gcd(n,k))/phi(n/gcd(n,k))*x^k:
od: RETURN(R[n](x)=a) end;
\end{verbatim}
\newpage

\centerline{Table of $R_n(x)$ for $1\le n\le 20$}
\[
\vbox {\offinterlineskip  \hrule \halign{ \strut \vrule
       \hfill $\ # \ $  & \vrule
      $\ #  $ \hfill  \vrule  \cr
      n & R_n(x) \ \cr \noalign{\hrule}
      1 &  1 \ \cr \noalign{\hrule}
      2 & 1-x \cr  \noalign{\hrule}
      3 & 2-x-x^{2}=(1-x)(2+x) \cr  \noalign{\hrule}
      4 & 2-2x^{2}=2(1-x)(1+x) \cr  \noalign{\hrule}
      5 & 4-x-x^{2}-x^{3}-x^{4}=(1-x)(4+3x+2x^2+x^3)\cr  \noalign{\hrule}
      6 & 2+x-x^{2} - 2 x^{3} - x^{4} + x^{5} \cr
        & = (1 - x)(2 - x)(1+x)(1+x+x^{2})\cr  \noalign{\hrule}
      7 & 6 - x - x^{2} - x^{3} - x^{4} - x^{5} - x^{6}\cr
        & = (1 -x)(6+5x+4x^2+3x^3+2x^4+ x^{5})\cr  \noalign{\hrule}
      8 & 4 - 4x^{4}= 4(1-x)(1+x)(1+x^{2})\cr  \noalign{\hrule}
      9 & 6 - 3x^{3} - 3x^{6}= 3(1-x)(2+x^{3})(1+x+x^{2})\cr  \noalign{\hrule}
     10 & 4 + x - x^{2} + x^{3} - x^{4} - 4x^{5} - x^{6} + x^{7} - x^{8} + x^{9} \cr
        & = (1-x)(1+x)(4-3x+2x^2- x^{3})(1+x+x^2+x^3+x^4) \cr \noalign{\hrule}
     11 & 10 - x - x^{2} - x^{3} - x^{4} - x^{5} - x^{6} - x^{7} - x^{8} - x^{9} -
           x^{10} \cr
        & = (1-x)(10 +9x +8x^2+7x^3 +6x^4+5x^5+4x^6+3x^7+2x^8+
           x^{9}) \cr \noalign{\hrule}
     12 & 4 + 2x^{2} - 2x^{4} - 4x^{6} - 2x^{8} +2x^{10} \cr
        &  = 2(1-x)(1+x)(2-x^{2})(1+x^{2})(1-x+x^{2})(1+x+x^{2}) \cr \noalign{\hrule}
     13 & 12 - x- x^{2} - x^{3} - x^{4} - x^{5} - x^{6} - x^{7} - x^{8} - x^{9} - x^{10} - x^{11} - x^{12} \cr
        & = (1-x)(12+11x+10x^2+9x^3+8x^4+7x^5+6x^6+5x^7 +4x^8 \cr
        & +3x^9+2x^{10} +x^{11})\cr \noalign{\hrule}
     14 & 6 + x - x^{2} + x^{3} - x^{4} + x^{5} - x^{6} - 6x^{7} - x^{8} + x^{9} - x^{10} + x^{11} - x^{12}
          + x^{13} \cr
        & = (1-x)(1+x)(1+x+x^2+x^3+x^4+x^5 + x^{6})(6- 5x + 4x^2-3x^3 \cr
        & +2x^4-x^{5})\cr \noalign{\hrule}
     15 & 8 + x + x^{2} - 2x^{3} + x^{4} - 4x^{5} - 2x^{6} + x^{7} + x^{8}
          - 2x^{9} - 4x^{10} + x^{11} - 2x^{12} + x^{13} + x^{14} \cr
        & = (1-x)(1+x+x^2+x^3+ x^{4})(8-7x+5x^3-4x^4+3x^5-x^{7})(1+x+x^{2}) \cr \noalign{\hrule}
     16 & 8 - 8x^{8} =  8(1-x)(1+x)(1+x^{2})(1+x^{4}) \cr \noalign{\hrule}
     17 & 16 - x - x^{2} - x ^{3} - x^{4} - x^{5} - x^{6} - x^{7} - x^{8} - x^{9} - x^{10} -
           x^{11} - x^{12} \cr
        & - x^{13} - x^{14} - x^{15} - x^{16} \cr
        & = (1-x)(16+15x+ 14x^{2} + 13x^{3} + 12x^{4} + 11x^{5} + 10x^{6} + 9x^{7} + 8x^{8} \cr & + 7x^{9} +
          6x^{10} + 5x^{11} + 4x^{12} + 3x^{13} + 2x^{14} + x^{15}) \cr \noalign{\hrule}
     18 & 6 + 3x^{3} - 3x^{6} - 6x^{9} - 3x^{12} + 3x^{15} \cr
        & = 3(1 -x)(1+x)(1-x+x^{2})(1+x+x^{2})(1+x^3+x^{6})(2-x^{3}) \cr \noalign{\hrule}
     19 & 18 - x - x^{2} - x^{3} - x^{4} - x^{5} - x^{6}
          - x^{7} - x^{8} - x^{9} - x^{10} - x^{11} - x^{12} \cr & - x^{13} - x^{14} - x^{15} - x^{16} - x^{17} -
           x^{18} \cr
        & =  (1-x)(18+17x+16x^2+15x^3+ 14x^{4} + 13x^{5} + 12x^{6} + 11x^{7} + 10x^{8} \cr
        & + 9x^{9} + 8x^{10} + 7x^{11} + 6x^{12} + 5x^{13} + 4x^{14} + 3x^{15} + 2x^{16} + x^{17})
           \cr \noalign{\hrule}
     20 & 8 + 2x^{2} - 2x^{4} + 2x^{6} - 2x^{8} - 8x^{10} - 2x^{12} + 2x^{14} - 2x^{16} + 2x^{18} \cr
        & = 2(1-x)(1+x)(1+x^{2})(4-3x^{2} + 2x^{4}-x^{6})(1+x+x^2+x^{3} +
          x^{4}) \cr
        & (1-x+x^2-x^3+x^{4})
\cr} \hrule}
\]

\newpage

\centerline{Table of $T_n(x)$ for $1\le n\le 20$}
\[
\vbox{\offinterlineskip  \hrule \halign{ \strut \vrule
       \hfill $\ # \ $  & \vrule
      $\ #  $ \hfill  \vrule  \cr
      n & T_n(x) \ \cr \noalign{\hrule}
      1 &  1 \ \cr \noalign{\hrule}
      2 & 1+x \cr  \noalign{\hrule}
      3 & 2+x+x^{2} \cr  \noalign{\hrule}
      4 & 2+2x^{2}=2(1+x^2) \cr  \noalign{\hrule}
      5 & 4+x+x^{2}+x^{3}+x^{4} \cr  \noalign{\hrule}
      6 & 2+x+x^{2}+2 x^{3}+ x^{4} + x^{5} \cr
        & = (1+x)(1-x+x^2)(2+x+x^{2}) \cr  \noalign{\hrule}
      7 & 6 + x + x^{2} + x^{3} + x^{4} + x^{5} + x^{6} \cr  \noalign{\hrule}
      8 & 4 + 4x^{4}= 4(1+x^{4}) \cr  \noalign{\hrule}
      9 & 6 + 3x^{3} + 3x^{6} \cr  \noalign{\hrule}
     10 & 4 + x + x^{2} + x^{3} + x^{4} + 4x^{5} + x^{6} + x^{7} + x^{8} + x^{9} \cr
        & = (1+x)(4+x+x^2+x^{3}+x^4)(1-x+x^2-x^3+x^4) \cr \noalign{\hrule}
     11 & 10 + x + x^{2} + x^{3} + x^{4} + x^{5} + x^{6} + x^{7} + x^{8} + x^{9}
           + x^{10} \cr \noalign{\hrule}
     12 & 4 + 2x^{2} + 2x^{4} + 4x^{6} + 2x^{8} + 2x^{10} \cr
        &  = 2(1+x^2)(2+x^{2}+x^4)(1-x^{2}+x^4) \cr \noalign{\hrule}
     13 & 12 + x + x^{2} + x^{3} + x^{4} + x^{5} + x^{6} + x^{7} + x^{8} + x^{9} + x^{10}+ x^{11}+ x^{12} \cr
          \noalign{\hrule}
     14 & 6 + x + x^{2} + x^{3} + x^{4} + x^{5} + x^{6} + 6x^{7} + x^{8} + x^{9} + x^{10} + x^{11} + x^{12}
          + x^{13} \cr
        & = (1+x)(1-x+x^2-x^3+x^4-x^5 + x^{6})(6+x + x^2+ x^3 +x^4+x^{5}+x^6)
        \cr \noalign{\hrule}
     15 & 8 + x + x^{2} + 2x^{3} + x^{4} + 4x^{5} + 2x^{6} + x^{7} + x^{8}
          + 2x^{9} + 4x^{10} + x^{11} + 2x^{12} + x^{13} + x^{14} \cr \noalign{\hrule}
     16 & 8 + 8x^{8} = 8(1+x^8) \cr \noalign{\hrule}
     17 & 16 + x + x^{2} + x ^{3} +x^{4} + x^{5} + x^{6} + x^{7} + x^{8} + x^{9} + x^{10}
          + x^{11} + x^{12} \cr
        & + x^{13} + x^{14} + x^{15} + x^{16} \cr \noalign{\hrule}
     18 & 6 + 3x^{3} + 3x^{6} + 6x^{9} + 3x^{12} + 3x^{15} \cr
        & = 3(1+x)(1-x+x^{2})(2+x^3+x^{6})(1-x^3+x^{6}) \cr \noalign{\hrule}
     19 & 18 + x + x^{2} + x^{3} + x^{4} + x^{5} + x^{6}
          + x^{7} + x^{8} + x^{9} + x^{10} + x^{11} + x^{12} \cr & + x^{13} + x^{14} + x^{15} + x^{16} + x^{17}
          + x^{18} \cr \noalign{\hrule}
     20 & 8 + 2x^{2} + 2x^{4} + 2x^{6} + 2x^{8} + 8x^{10} + 2x^{12} + 2x^{14} + 2x^{16} + 2x^{18} \cr
        & = 2(1+x^2)(4+x^{2} + x^{4}+x^{6}+x^8)(1-x^2+x^4-x^{6} +x^{8})
\cr} \hrule}
\]

\end{document}